\documentclass[a4paper]{amsart}
\usepackage{amssymb}
\usepackage{amsmath}
\usepackage{graphicx}
\usepackage{amsthm}
\usepackage{dsfont}
\usepackage{color}
\begin{document}
\thanks{The first author was
 partially supported by  PSC CUNY Research Award,
 No. 60007-33-34, and NSF grant  DMS 0401318}
\newtheorem{introtheorem}{Theorem}
\renewcommand{\theintrotheorem}{\Alph{introtheorem}}
\newtheorem{theorem}{Theorem }[section]
\newtheorem{lemma}[theorem]{Lemma}
\newtheorem{corollary}[theorem]{Corollary}
\newtheorem{proposition}[theorem]{Proposition}
\theoremstyle{definition}
\newtheorem{definition}[theorem]{Definition}
\newtheorem{example}[theorem]{Example}
\theoremstyle{remark}
\newtheorem{remark}[theorem]{Remark}

\renewcommand{\labelenumi}{(\roman{enumi})} 
\def\theenumi{\roman{enumi}}

\numberwithin{equation}{section}

\def \g {{\gamma}}
\def \G {{\Gamma}}
\def \l {{\lambda}}
\def \a {{\alpha}}
\def \b {{\beta}}
\def \f {{\phi}}
\def \r {{\rho}}
\def \R {{\mathbb R}}
\def \H {{\mathbb H}}
\def \N {{\mathbb N}}
\def \C {{\mathbb C}}
\def \Z {{\mathbb Z}}
\def \F {{\Phi}}
\def \Q {{\mathbb Q}}
\def \e {{\epsilon }}
\def \ev {{\vec\epsilon}}
\def \ov {{\vec{0}}}
\def \GinfmodG {{\Gamma_{\!\!\infty}\!\!\setminus\!\Gamma}}
\def \GmodH {{\Gamma\backslash\H}}
\def \sl  {{\hbox{SL}_2( {\mathbb R})} }
\def \psl  {{\hbox{PSL}_2( {\mathbb R})} }
\def \L  {{\hbox{L}^2}}

\newcommand{\norm}[1]{\left\lVert #1 \right\rVert}
\newcommand{\abs}[1]{\left\lvert #1 \right\rvert}
\newcommand{\modsym}[2]{\left \langle #1,#2 \right\rangle}
\newcommand{\inprod}[2]{\left \langle #1,#2 \right\rangle}
\newcommand{\Nz}[1]{\left\lVert #1 \right\rVert_z}
\newcommand{\ver}[1]{\operatorname{vert}\left( #1 \right)}
\newcommand{\wl}[1]{\operatorname{wl}\left( #1 \right)}
\newcommand{\wlc}[1]{\operatorname{wl_c}\left( #1 \right)}

\title{Discrete logarithms in free groups}
\author{Yiannis N. Petridis}
\address{Department of Mathematics and Computer Science\\
City University of New York, Lehman College\\
250 Bedford Park Boulevard West
Bronx\\NY 10468-1589\newline\mbox{} \hspace{8pt} The Graduate Center, Mathematics Ph.D. Program\\
                       365 Fifth Avenue, Room 4208
                       New York, NY 10016-4309}
\email{petridis@comet.lehman.cuny.edu}
\author{Morten S. Risager}
\address{Department of Mathematical Sciences\\ University of Aarhus\\ Ny Munkegade Building 530\\ 8000 {Aa}rhus, Denmark}
\email{risager@imf.au.dk}
\date{\today}

\begin{abstract}For the free group on $n$ generators we prove that the
  discrete logarithm is distributed according to the standard Gaussian
  when the logarithm is renormalized appropriately. 
\end{abstract}
\maketitle
\section{Introduction}  
Let $\G=F(A_1,\ldots,A_n)$, $n\geq 2$ be the free group on $n$ generators. We define additive homomorphisms on the generators
\begin{equation}\begin{array}{rccc}
\log_j: &\G&\to &\Z\\
        &A_i&\mapsto&\delta_{ij}.
\end{array}\end{equation}
Hence $\log_j$ counts the number of occurrences (with signs) of the
generator $A_j$. We want to study the distribution of this discrete
logarithm as the
cyclically reduced word length grows to infinity.  Our main result is
that after a suitable re-normalization and restriction the discrete logarithm is
distributed according to a standard Gaussian distribution. More precisely, let 
\begin{equation} 
\mathfrak{log}_j(\g)=\sqrt{\frac{n-1}{\wl{\g}}}\log_j(\g).
\end{equation} Here $\wl{\g}$ denotes the word
length of $\g\in\G$. Let $\G_c$ be the set of cyclically reduced words in $\G$. 
\begin{theorem}\label{mainresult} In $\G_c$ the function $\mathfrak{log}_j$ has
  asymptotically a standard Gaussian distribution. More precisely, we have 
\begin{equation*}
\frac{\#\{\g\in\G_c|\wl{\g}\leq l, \,
  \mathfrak{log}_j(\g)\in [a,b]\}}{\#\{\g\in\G_c|\wl{\g}\leq l\}}\to\frac{1}{\sqrt{2\pi}}\int_a^b\exp\
    \left(-\frac{x^2}{2}\right) dx
  \textrm{ as }l\to\infty
\end{equation*}
\end{theorem}
After proving this theorem we discovered that I. Rivin arrived at a
similar theorem in \cite[Theorem 5.1]{rivin} (see also \cite{sharp}).
Rivin's proof  uses certain results on Chebyshev polynomials. Our proof uses character perturbations of the adjacency operator of a
singleton graph (here when $n=4$): 
\begin{figure}[!ht]
\centering
\includegraphics{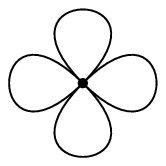}
%\caption{n=4}
\end{figure} 
\\ There is a fundamental identity, due to Ihara, which relates geometric
entities (lengths of closed paths on the graph) to the spectrum of
the adjacency operator of the graph. We consider this identity when
the system transforms according to a certain smooth family of characters
related to the  discrete logarithms. By differentiating in the character
family  we can extract enough information about
the relevant generating functions to calculate the moments of the
random variable in Theorem \ref{mainresult}. Apart from the use of this identity our proof is elementary. Our method is 
very powerful and has been used by the authors in a series of papers
\cite{petridis,petr, risager, risagerpetr} to prove distribution
results of additive homomorphisms in various different contexts. In fact we were motivated by our previous results to investigate the normal distribution on the level  of the group, without considering its action as isometries on hyperbolic space.
\section{From free groups to graphs}\label{fromto}
Let $\G=F(A_1,\ldots,A_n)$ $n\geq 2$ be the free group on $n$
generators. We want to relate a finite graph to $\G$. We refer to
\cite{serre,terras, nikitinvenkov,nikitin} for the basics on graphs and further explanation
of some of the terminology that we use. We mostly follow the
terminology of \cite{nikitin}.

Every $g\in\G$ is uniquely represented by a reduced word i.e.
$$\g=\prod_{i=1}^k A_{m_i}^{\epsilon_i}\qquad m_i\in\{1,\ldots,n\}, \epsilon_i\in\{\pm 1\} $$
in which there does not exist $i_0$ such that $m_{i_0}= m_{i_0+1}$ and
$\epsilon_{i_0}\epsilon_{i_0+1}=-1$. Hence 
\begin{equation}\log_j(\g)=\sum_{\substack{i=1\\m_i=j}}^k\e_{m_i}
\end{equation}
counts (with signs) the number of occurrences of $A_j$ in the reduced word representation of $\g$. We define the word length of $\g$ to be $\wl{\g}=k$. We also define $g\in\G$ to be cyclically reduced if $$\g\circ\gamma=\prod_{i=1}^k A_{m_i}^{\epsilon_i}\prod_{i=1}^k A_{m_i}^{\epsilon_i}$$ is reduced (or equivalently if either $\e_1\e_k\neq -1$ or $m_1\neq m_k$). We consider the set $\G_c$ of cyclically reduced words $$\G_c=\{\g\in\G\vert \g \textrm{ is cyclically reduced}\}.$$
There is an infinite-to-one map $\psi$ (cyclic reduction) from $\G$ to $\G_c$ and this map satisfies  $\log_{j}(\g)=\log_{j}(\psi(\g)).$

Let $S=\{A_1,\ldots,A_n\}$ and consider the $2n$-regular Cayley graph $X=X(\G,S)$ i.e. the graph having vertices $\G$ and (positive) edges $\G\times S$ where the edge $(g,s)\in \G\times S$ has the following origin and terminus \begin{eqnarray*}o(g,s)&=&g\\t(g,s)&=&gs.\end{eqnarray*} The group $\G$ acts (strictly hyperbolic) by isometric isomorphisms on $X$ and the quotient graph $\G\backslash X$ is $2n$-regular and has one vertex.

We denote by $ \operatorname{C}^{\operatorname{red}}$ the set of
closed reduced paths on the graph. The set of elements in  
$\operatorname{C}^{\operatorname{red}}$ of length $n$ is denoted by 
$\operatorname{C}^{\operatorname{red}}_n$. For this particular graph
the following is easy
to verify:
\begin{lemma}\label{basiccorrespondence}There is a one-to-one
  correspondence between $\operatorname{C}^{\operatorname{red}}_m$ and
  the set of cyclically reduced group elements $\g\in \G_c$ with $\wl{\g}=m$.\end{lemma} 
We can now therefore define an additive homomorphism $\widetilde\log_j$ on
  $\operatorname{C}^{\operatorname{red}}$ by fixing an isomorphism
  $\phi: \operatorname{C}^{\operatorname{red}}_1\to \{\g\in \G_c| \wl{\g}=1\}$ and defining
  $\widetilde\log_j(C)=\log_j(\phi(C))$. By an obvious abuse of notation
  we shall often write $\log_j(C)$ instead of $\widetilde\log_j(C)$.

\section{The Ihara zeta function}
We now briefly explain the basics of the Ihara zeta function. There
is a close resemblance between the Ihara zeta function and the Selberg
zeta function \cite{hejhal}. We
refer to \cite{terras,nikitinvenkov, nikitin} for further details.
Let $X$ be a $(q+1)$-regular infinite tree ($q\geq 2$) and let $\G$ be a strictly hyperbolic free subgroup of isometric automorphisms of $X$ with finite $(q+1)$-regular quotient graph \begin{equation}\#\ver{\G\backslash X}< \infty.\end{equation} Let furthermore $\chi:\Gamma\to S^1$ be a unitary character. Consider
\begin{equation*}L^2(X,\chi)=\left\{f:\ver{X}\to\C\left\vert  f(\g v)=\chi(\g)f(v)\textrm{ for all }\g\in\G, v\in\ver{X}\right.\right\}
\end{equation*}
together with the inner product \begin{equation*}\inprod{f}{g}=\sum_{v\in\ver{\G\backslash X}}f(v)\overline{g(v)}.\end{equation*} We note that this is a finite-dimensional Hilbert space of dimension $\dim{L^2(X,\chi)}=\#\ver{\G\backslash X}$. We consider the adjacency operator $A(\G,\chi):L^2(X,\chi)\to L^2(X,\chi)$ defined by \begin{equation}(A(\G,\chi)f)(v)=\sum_{v'\sim v}f(v').\end{equation} Here $v'\sim v$ means that $v'$ is adjacent to $v$. The adjacency operator is self-adjoint with spectrum contained in $[-(q+1),q+1]$. 

Let $\mathfrak{P}_\G$ be the set of primitive conjugacy classes in $\G$ different from the unit class. For $\{P\}\in\mathfrak{P}_\G$ we define the degree \begin{equation}\deg(P)=\min_{v\in\ver{X}}l(v,Pv)\end{equation} where $l(v,v')$ is the length of the geodesic from $v$ to $v'$.

The Ihara zeta function is defined, for $u\in \C$ with $\abs{u}<q^{-1}$, by 
\begin{equation}Z(\G,\chi; u)=\prod_{\{P\}\in\mathfrak{P}_\G}(1-\chi(P)u^{\deg{P}})^{-1}.\end{equation} There is a fundamental identity, due to Ihara, which relates this geometric object to the spectrum of the adjacency operator. More precisely, we have the following explicit determinant representation of the Ihara zeta function \cite{nikitinvenkov}:
\begin{theorem}[Ihara]
\begin{equation}\label{iharaidentity}
Z(\G,\chi; u)=(1-u^2)^{-g}\det{(1-uA(\G,\chi)+qu^2)}^{-1},
\end{equation}
where $g=(q-1) \# \ver{\G\backslash X}/2$. 
\end{theorem}
This gives immediately that the Ihara zeta function admits meromorphic continuation to the whole complex plane.
We consider the logarithmic derivative of $Z(\G,\chi; u)$ and find
\begin{eqnarray}u\frac{d}{d u}\log Z(\G,\chi;u)&=& u \sum_{\{P\}\in\mathfrak{P}_\G}\frac{\chi(P)\deg{P} u^{\deg{P}-1}}{1-\chi(P)u^{\deg{P}}}\\
\nonumber &=&\sum_{n=1}^\infty\sum_{\{P\}\in\mathfrak{P}_\G}\chi(P^n)\deg{P}u^{n\deg{P}}\\\nonumber &=&\sum_{m=1}^\infty  n_{\G,\chi}(m)u^m,
\end{eqnarray} 
where
\begin{equation}
n_{\G,\chi}(m)=\sum_{\substack{n\deg{P}=m\\ n\in\N, \{P\}\in\mathfrak{P}_\G}}\chi(P^n)\deg{P}.
\end{equation}
There is a one-to-one correspondence between primitive conjugacy classes of $\G$ of degree $d$, $\mathfrak{P}_\G(d)$ and primitive cycles in $\G\backslash X$ of  length $d$, $\operatorname{Prim}_{\G\backslash X}(d)$ . The natural map from the set of closed reduced prime paths of length $d$ in $\G\backslash X$, $\operatorname{C}^{\operatorname{red},\operatorname{prime}}_d$, to the set of primitive cycles of length $d$ in $\G\backslash X$ is $d$-to-one (each of the $d$ vertices can be a starting point for the cycle). We therefore have
\begin{equation*}n_{\G,\chi}(m)=\sum_{d|m}d\sum_{\{P\}\in\mathfrak{P}_{\G}}\chi(P^{m/d})=\sum_{d|m}\sum_{C\in\operatorname{C}^{\operatorname{red},\operatorname{prime}}_d}\chi(C^{m/d})= \sum_{C \in \operatorname{C}^{\operatorname{red}}_m}\chi(C).
\end{equation*}
Using (\ref{iharaidentity}) we get 
\begin{equation}\label{niceidentity} \sum_{m=1}^\infty n_{\G,\chi}(m)u^m=\frac{2gu^2}{1-u^2}-u\frac{\frac{d}{du}\det(1-uA(\G,\chi)+qu^2)}{\det(1-uA(\G,\chi)+qu^2)}
\end{equation}
When $\chi=1$ the adjacency operator has $q+1$ as an eigenvalue with
the  identity as eigenfunction. Hence, in this case, we see that (\ref{niceidentity}) has a pole at $u=q^{-1}$.

\section{Generating functions}
We now explain how to  use the Ihara zeta function to get
information about the generating functions relevant to Theorem \ref{mainresult}. Since we are ultimately interested in the case of the
singleton graph constructed in section \ref{fromto} we assume for
simplicity that the quotient in the previous section has only one vertex. 
In the case of the graph from section \ref{fromto} we have $q+1=2n$.
 
We define a family of unitary characters on $\G$ by $$\chi_\e(\g)=\exp(i\e\log_j(\g)).$$ In this simple situation the adjacency operator is simply a multiplication operator:
\begin{equation}(A(\G,\chi_\e)f)(v)=\sum_{i=1}^n(\chi_\e(A_i)+\chi_\e(A_i)^{-1})f(v).\end{equation} For simplicity we let \begin{equation}A(\e)=\sum_{i=1}^n(\chi_\e(A_i)+\chi_\e(A_i)^{-1})=2(n-1)+2\cos (\e ).\end{equation} Equation (\ref{niceidentity}) now takes the form
\begin{equation}\label{willihaveajob}
 \sum_{m=1}^\infty n_{\G,\chi_\e}(m)u^m=\frac{2gu^2}{1-u^2}+\frac{uA(\e)-2qu^2}{1-uA(\e)+qu^2}
\end{equation}
We denote the left-hand-side by $G(u,\e)$. This is the generating function relevant to the problem we are studying. When $\e=0$ this function has a simple pole at $u=q^{-1}$ with residue 
\begin{equation}\label{originalresidue}\operatorname{res}_{u=q^{-1}}G(u,0)=-q^{-1}.\end{equation}
We wish to understand
\begin{equation}G^{(k)}(u):=\left.\frac{d^k}{d\e^k} G(u,\e)\right\vert_{\e=0}=\sum_{m=1}^\infty\left(\sum_{C \in \operatorname{C}^{\operatorname{red}}_m}i^k\log_j(C)^k\right)u^m 
\end{equation}
The convergence is uniform in $\e$, by comparison with the series with $\e=0$, so we can differentiate termwise.
%\begin{center} *Need to say why we can differentiate termwise***\end{center}

From (\ref{willihaveajob}) we get immediately the following result:
\begin{theorem}\label{continuation} The function $G^{(k)}(u)$ admits
  meromorphic continuation to the whole complex plane. The points
  $u=q^{-1}$, $u=1$ are the only possible poles when $n\geq 1$. The function $G(u)=G^{(0)}(u)$ has poles at $u=\pm 1$ and $u=q^{-1}$.\end{theorem}

We now calculate the pole order and leading term of the pole at $u=q^{-1}$.
When $k\ge 1$ Eq. (\ref{willihaveajob})  gives 
\begin{equation}\label{bluesky}
G^{(k)}(u)=\sum_{l=0}^k\binom{k}{l}(uA^{(l)}(0)-\delta_{l=0}2qu^2)\left.\frac{d^{k-l}}{d\e^{k-l}}(1-uA(\e)+qu^2)^{-1}\right\vert_{\e=0}
\end{equation}
Clearly we have
\begin{equation}\label{Acalculated}
A^{(l)}(0)=\begin{cases}2n,&\textrm{ if $l=0$,}\\0, &\textrm{ if $l$ is odd,}\\2(-1)^{l/2},& \textrm{ if $l\geq 2$ is even. }\end{cases}
\end{equation}
Since
 \begin{equation}(1-uA(\e)+qu^2)(1-uA(\e)+qu^2)^{-1}=1, \end{equation}
 we find that 
\begin{align}\label{pencil}&\left.\frac{d^{l}}{d\e^l}(1-uA(\e)+qu^2)^{-1}\right\vert_{\e=0}\\&= (1-uA(0)+qu^2)^{-1}\sum_{r=0}^{l-1}\binom{l}{r}A^{(l-r)}(0)u\left.\frac{d^{r}}{d\e^r}(1-uA(\e)+qu^2)^{-1}\right\vert_{\e=0}\nonumber\end{align}
When $l=1$ it is clear from (\ref{Acalculated}) that this vanishes. By induction and (\ref{Acalculated})  we conclude that 
\begin{equation}\label{pen}
\left.\frac{d^{l}}{d\e^l}(1-uA(\e)+qu^2)^{-1}\right\vert_{\e=0}=0
\end{equation} when $l$ is odd. 
Alternatively $1-uA(\e )+qu^2$ is an even function.

From (\ref{pencil}), (\ref{pen}), and (\ref{Acalculated}) an easy induction argument now shows that when $l$ is even  the left-hand-side of  (\ref{pencil}) has a pole of order $l/2+1$ and the leading term, $L_l$, satisfies the recurrence relation
\begin{equation}
L_l=(\operatorname{res}_{u=q^{-1}}(1-A(0)u+qu^2)^{-1})\binom{l}{l-2}A^{(2)}(0)q^{-1}L_{l-2}.
\end{equation}  From this we find easily
\begin{equation}
L_l=-\frac{1}{q^{l/2}}\frac{l!}{(q-1)^{l/2+1}}.
\end{equation}
We summarize our  result:
\begin{theorem}\label{laurent} When $k$ is odd the function $G^{(k)}(u)$ is identically zero. When $k$ is even $G^{(k)}(u)$ has a pole at $u=q^{-1}$ of order $k/2+1$ and the leading term in the Laurent expansion around $u=q^{-1}$ equals
$$-\frac{k!}{q^{k/2+1}(q-1)^{k/2}}.$$
\end{theorem}
We notice that the vanishing of $G^{(k)}(u)$ for $k$ odd is also clear from the fact that $$\sum_{C \in \operatorname{C}^{\operatorname{red}}_m}\log_j(C)^k=0.$$
This fact is  clear since if $C$ is a reduced path of length $m$ then the inverse path is also a reduced path of length $m$ and $\log_j(C)^k+\log_j(C^{-1})^k=0$ when $m$ is odd. 
\section{Summations}
In this section we want to find asymptotics as $T\to\infty$ of the following sums
\begin{equation}\label{wewantthis}
\sum_{\substack{\g\in\G_c \\ \wl{\g}\leq T}}\log_j(\g)^k
\end{equation} 
In principle we can find exact expressions for these sums. We show how
this is done when $k=0$. In the general situation we settle with
asymptotics since this is all we need  for Theorem \ref{mainresult}. We do this
on the basis of the calculations of the previous chapter and a
Tauberian theorem to be proved below.
We notice that from Lemma \ref{basiccorrespondence} and the fact that if $\phi(C)=\g$ then $l(C)=\wl\g$ we have that (\ref{wewantthis}) equals
\begin{equation}\label{alternative}\sum_{m\leq T}\sum_{C\in \operatorname{C}^{\operatorname{red}}_m}\log_j(C)^k\end{equation} For $k$ odd the discussion in the end of the last section shows that (\ref{wewantthis}) is identically zero. 

We have 
$$\frac{1}{1-u2n+qu^2}=\frac{1}{q-1}\left(\frac{q}{1-qu}-\frac{1}{1-u}\right).$$
From this we see, using a geometric expansion, that 
$$\frac{u2n-2qu^2}{1-u2n+qu^2}=\sum_{m=1}^{\infty}(q^m+1)u^m.$$
Since
$$\frac{2gu^2}{1-u^2}=\frac{q-1}{2}\sum_{m=1}^\infty\left(1+(-1)^m\right)u^m,$$
we conclude that
$$G^{(0)}(u)=\sum_{m=1}^\infty \#
\operatorname{C}^{\operatorname{red}}_m u^m
=\sum_{m=1}^\infty\left(q^m+1+\frac{q-1}{2}\left(1+(-1)^m\right)\right)u^m
$$
when $\abs{u}<q^{-1}$. Using Lemma \ref{basiccorrespondence}  and the uniqueness of Taylor expansions, we proved the following theorem:
\begin{theorem}Let $\G$ be the free group on $n$ elements and $\G_c$ the set of cyclically reduced words in $\G$. Then the
  number of cyclically reduced words of word length $m$ equals
$$\#\{\g\in\G_c| \wl{\g}=m\}=(2n-1)^m+1+(n-1)\left(1+(-1)^m\right)$$
\end{theorem}
We notice that this is \cite[Theorem 1.1]{rivin}. In principle we
can use the same technique as above  for any $k$ using
(\ref{bluesky}). This would give explicit - though quite complicated - expressions for  
$$\sum_{C\in
  \operatorname{C}^{\operatorname{red}}_m}\log_j(C)^k.$$ For
  simplicity we use only the fact that we have calculated the leading
  term of (\ref{bluesky}) to obtain main term and sharp error
  estimates for (\ref{alternative})

To do this we need the following: Let 
$$f(u)=\sum_{m=0}^\infty c_m u^m$$ have radius of convergence
$R= q^{-1}\leq 1$. Assume that $f$ admits meromorphic continuation to an open set
containing $\{u\in\C| \abs{u}\leq 1\}$. We assume that $u=q^{-1}$ and
$u=1$ are the only possible poles in $\{u\in\C| \abs{u}\leq 1\}$. Let
\begin{equation*}
\sum_{k=1}^{k_1}\frac{a_k}{(u-q^{-1})^k}\quad \textrm { resp. }\quad \sum_{k=1}^{k_2}\frac{b_k}{(u-1)^k}
\end{equation*} be the singular parts of $f$ at $u=q^{-1}$
resp. $u=1$. Now the function 
\begin{equation}f(u)-\sum_{k=1}^{k_1}\frac{a_k}{(u-q^{-1})^k}-\quad
  \sum_{k=1}^{k_2}\frac{b_k}{(u-1)^k}
\end{equation} is regular in an open disc containing the closure of the unit disc. 
Using the series expansions
\begin{eqnarray*}\frac{1}{(u-1)^k}&=&(-1)^k\sum_{m=0}^\infty\frac{(k)_m}{m!}u^m,\\
\frac{1}{(u-q^{-1})^k}&=&q^k(-1)^k\sum_{m=0}^\infty\frac{(k)_mq^m}{m!}u^m,
\end{eqnarray*}
where $(k)_m=k(k+1)\cdots(k+(m-1))$ we conclude that the series
\begin{equation*}\sum_{m=0}^\infty \left(c_m -\sum_{k=1}^{k_1}a_k(-q)^k\frac{(k)_m}{m!}q^m-\sum_{k=1}^{k_2}b_k(-1)^k\frac{(k)_m}{m!} \right)u^m\end{equation*}
has radius of convergence strictly larger than 1. In particular it is
convergent at $u=1$, that is 
\begin{equation}\label{ugly}\sum_{m=0}^l c_m
  =\sum_{m=0}^l\sum_{k=1}^{k_1}a_k(-q)^k\frac{(k)_m}{m!}q^m+\sum_{m=0}^l\sum_{k=1}^{k_2}b_k(-1)^k\frac{(k)_m}{m!}
  +o(1) \end{equation} as $l\to\infty$.
We then use the following elementary lemma:
\begin{lemma} Let $q>1$ and $k\in \N$. We have the following asymptotic expansions
\begin{align}
\label{firstone}\sum_{m=0}^l\frac{(k)_m}{m!}=&\frac{1}{k!}l^{k}+O(l^{k-1})\\
\label{secondone}\sum_{m=0}^l\frac{(k)_m}{m!}q^m=&\frac{1}{(q-1)(k-1)!}q^{l+1}l^{k-1}+\begin{cases}O(q^{l+1}l^{k-2}),
&  \textrm{ if }k\geq 2,\\O(1) &\textrm{ if }k=1.\end{cases}
\end{align}
\end{lemma}
\begin{proof}
To prove (\ref{firstone}) we notice that
$(k)_m/m!=(m+1)(m+2)\cdots(m+k-1)/(k-1)!$. Using the elementary estimate 
\begin{equation*}
\sum_{m=0}^{l}m^{k-1}=\frac{1}{k}l^k+O(l^{k-1})
\end{equation*}
proves the claim.

To prove (\ref{secondone}) we define $$a(k,m,q)=\frac{(k)_m}{m!}q^m \quad\textrm{ and } \quad
A(k,l,q)=\sum_{m=0}^la(k,m,q).$$ Then using 
$a(k,m,q)=(k+m-1)a(k-1,m,q)/(k-1)$ (when $k\neq 1$) and partial summation we find
\begin{align}\label{induction}
\nonumber A(k,l,q)&=\frac{1}{k-1}\sum_{m=0}^la(k-1,m,q)(k+m-1)\\
&=\frac{1}{k-1}A(k-1,l,q)(k+l-1)-\frac{1}{k-1}\int_{0}^l A(k-1,t,q)dt\\
\nonumber &=\frac{l}{k-1}A(k-1,l,q)+A(k-1,l,q)-\frac{1}{k-1}\int_{0}^l A(k-1,t,q)dt
\end{align}
 When $k=1$ the claim is obvious since in this case
 $(k)_m/m!=1$. The remaining cases now follow from (\ref{induction})
 by induction.
\end{proof}

From this lemma and (\ref{ugly}) now follows that 
\begin{equation}\label{nice}
\sum_{m=0}^l c_m=\frac{(-q)^{k_1}a_{k_1}}{(q-1)(k_1-1)!}q^{l+1}l^{k_1-1}+O(q^{l+1}l^{k_1-2})
\end{equation}

We notice that we have proved a Wiener-Ikehara-type Tauberian theorem (See \cite{korevaar}) for
power series with \emph{general} coefficients under the assumption that the
sum admits meromorphic continuation to a disk containing the closure
of the unit disc.    

Using this, Theorem \ref{continuation}, and Theorem \ref{laurent} we
get 
\begin{theorem}\label{asymptotics}Let
\begin{equation}
\sum_{m=0}^l\sum_{C \in
  \operatorname{C}^{\operatorname{red}}_m}\log_j(C)^k=\begin{cases}0&\textrm{
  if }k\textrm{ is odd.}\\ \frac{k!}{(k/2)!}\frac{1}{(q-1)^{k/2+1}}q^{l+1}l^{k/2}+O(q^{l+1}l^{k/2-1})&\textrm{
  if }k\textrm{ is even.}\end{cases}
\end{equation}
\end{theorem}

We now define the normalized discrete logarithms to be
$$\mathfrak{log}_j(C)=\sqrt{\frac{g}{l(C)}}\log_j(C).$$
Here $l(C)$ is the length of the cycle $C$. We then define the random
variable $X_T$ with probability measure
\begin{equation*}
P(X_T\in [a,b])=\frac{\#\{C\in
  \operatorname{C}^{\operatorname{red}}|l(C)\leq l, \,
  \mathfrak{log}_j(C)\in [a,b]\}}{\#\{C\in
  \operatorname{C}^{\operatorname{red}}|l(C)\leq l\}}
\end{equation*}
We then calculate the asymptotical moments of these random variable
i.e. we find the asymptotics of 
$$M_k(X_l)=\frac{1}{\#\{C\in
  \operatorname{C}^{\operatorname{red}}|l(C)\leq l\}}\sum_{\substack{C\in
  \operatorname{C}^{\operatorname{red}}\\l(C)\leq l}}\mathfrak{log}_j^k$$
From Theorem \ref{asymptotics} we find, using partial summation, that
  as $l\to\infty$ 
\begin{equation}M_k(X_l)\to\begin{cases}\frac{k!}{2^{k/2}(k/2)!}&
  \textrm{if }k\textrm{ is even}\\
0& \textrm{if }k\textrm{ is even.}\end{cases}
\end{equation}
The left-hand side equals the moments of the standard Gaussian and from
a classical result due to  Fr\'echet and Shohat (See \cite[11.4.C]{loeve}) we conclude that 
\begin{equation}P(X_l\in [a,b])\to\frac{1}{\sqrt{2\pi}}\int_a^b\exp\
    \left(-\frac{x^2}{2}\right) dx
  \textrm{ as }l\to\infty \end{equation}
Using the identification in Lemma \ref{basiccorrespondence} and the fact that if $\phi(C)=\g$ then $l(C)=\wl\g$ we arrive
  at Theorem \ref{mainresult}.

{\bf Acknowledgments:}\newline{We would like to thank  Alexei B. Venkov and S\o ren Galatius for  valuable comments and 
suggestions. }


\begin{thebibliography}{QWEWE}


%\bibitem{ledrappier}  M. Babillot, F. Ledrappier: Lalley's Theorem on periodic orbits of hyperbolic flows.
%Ergodic Theory Dynam Systems,    {\bf 18} 1998, no. 1, 17--39.  





%\bibitem{galambos}J. Galambos, Distribution of arithmetical functions. A survey.  Ann. Inst. H. Poincar\'e Sect. B (N.S.) {\bf 6} (1970), 281--305;



%\bibitem{goldfeld1}  D. Goldfeld, Zeta functions formed with modular symbols. 
%Automorphic forms, automorphic representations, and arithmetic (Fort Worth, TX
% 1996), 111--121, Proc. Sympos. Pure Math., 66, Part 1, Amer.
%Math. Soc., Providence, RI, 1999. 


%\bibitem{goldfeld2} D.  Goldfeld, 
%The distribution of modular symbols. 
%Number theory in progress, Vol. 2 (Zakopane-Ko\'scielisko, 1997), 849--865, 
%de Gruyter, Berlin, 1999. 


%\bibitem{goldstein} L. J. Goldstein, Dedekind sums for s Fuchsian group. I.
%Nagoya Math. J. {\bf 50} (1973), 21--47.

%\bibitem{gradstein} I. D. Gradsh teyn, I. M. Ryzhik, Tables of Integrals, 
%Series, and Products, 5th ed., ed. Alan Jeffrey, Academic Press, Boston, 1994.

\bibitem{hejhal} D. Hejhal, The Selberg trace formula for ${\rm PSL}(2,\,R)$. 
Vol. 1. Lecture Notes in Mathematics, 1001.
Springer-Verlag, Berlin, 1983. viii+806pp.



%\bibitem{iwaniec}H. Iwaniec, {\it Introduction to the spectral theory of automorphic forms}, Rev. Mat. Iberoam., Madrid, 1995

\bibitem{korevaar} J. Korevaar,  A century of complex Tauberian theory. Bull. Amer. Math. Soc. (N.S.) {\bf 39} (2002), no.~4, 475--531 (electronic). 



%\bibitem{lalley}
%S. Lalley:
%Closed geodesics in homology classes on surfaces of variable negative curvature.
%Duke Math. J. 58 (1989), no. 3, 795--821.


\bibitem{nikitin}A. M. Nikitin, The Ihara-Selberg zeta function of a finite graph and symbolic dynamics, Algebra i Analiz {\bf 13} (2001), no. 5, 134--149; translation in St. Petersburg Math. J. {\bf 13} (2002), no.~5, 809--820.

\bibitem{nikitinvenkov}A. M. Nikitin and A. B. Venkov, The Selberg trace formula, Ramanujan graphs and some problems in mathematical physics. (Russian) , Algebra i Analiz {\bf 5} (1993), no. 3, 1--76; translation in St. Petersburg Math. J. {\bf 5} (1994), no.~3, 419--484. 




\bibitem{loeve} M. Lo\`eve, {\it Probability theory. I}, Fourth edition, Springer, New York, 1977.



\bibitem{petridis} Y. N. Petridis, Spectral deformations and Eisenstein Series Associated with Modular Symbols. Internat. Math. Res. Notices {\bf 2002}, 
no. 19, 991--1006.

\bibitem{petr} Y. N. Petridis, M. S. Risager, Modular symbols have a
  normal distribution, To appear in  Geom. Funct. Anal.

\bibitem{risagerpetr} Y. N. Petridis, M. S. Risager, The distribution of
  values of the Poincar\'e pairing for hyperbolic Riemann surfaces, To
  appear in  J. Reine ang. Mat.


%\bibitem{phillipssarnak1}
%R.  Phillips, P.  Sarnak, Cusp forms for character varieties. Geom. Funct. 
%Anal. 4 (1994), no. 1, 93--118.

%\bibitem{phillipssarnak2}
%R.  Phillips, P. Sarnak, The spectrum of Fermat curves. Geom. Funct. Anal. 1 
%(1991), no. 1, 80--146. 
 
%\bibitem{phillips3} R. Phillips, P. Sarnak, Geodesics in homology classes. Duke Math. J. {\bf 55} (1987), no. 2, 287--297. 

\bibitem{risager} M. S. Risager, On the distribution of modular symbols for
  compact surfaces, To appear in  Internat. Math. Res. Notices.


\bibitem{rivin}I. Rivin, Growth in free groups (and other stories), arXiv:math.CO/9911076. 
 
\bibitem{serre}J.-P. Serre, {\it Trees}, Translated from the French original by John Stillwell, Corrected 2nd printing of the 1980 English translation, Springer, Berlin, 2003.

\bibitem{sharp}R. Sharp, Local limit theorems for free groups, J
  Math. Ann, {\bf 321 }, 2001, 4, p. 889--904.

\bibitem{terras}A. Terras, {\it Fourier analysis on finite groups and applications}, Cambridge Univ. Press, Cambridge, 1999.

%\bibitem{selberg2}
% A. Selberg, Harmonic Analysis. G\"ottingen Lecture notes, in {\it Collected papers. Vol. I}, Springer, Berlin, 1989. 

%\bibitem{sharp} R. Sharp, 
% Closed geodesics and periods of automorphic forms. 
%Adv. Math. {\bf 160} (2001), no. 2, 205--216. 

%\bibitem{sharp2} R. Sharp, A local limit theorem for closed geodesics and homology, to appear in Trans. Amer. Math. Soc.

%\bibitem{venkov} A. Venkov,
%Spectral theory of automorphic functions. A translation of Trudy Mat. Inst. Steklov. 153 (1981). Proc. Steklov Inst. Math. 1982, no. 4(153), ix+163 pp. 1983.

%\bibitem{venkov2} A.  Venkov, {\it Spectral theory of automorphic functions and its applications}, Translated from the Russian by N. B. Lebedinskaya, Kluwer Acad. Publ., Dordrecht, 1990.


\end{thebibliography}
\end{document}